\title{ ~~\\ Value distribution of cyclotomic polynomial coefficients}
\author{Yves Gallot, Pieter Moree and Huib Hommersom}
\documentclass[12pt]{article}
\usepackage{amssymb, latexsym, amsfonts, makeidx}
\def\@ptsize{2}
\newtheorem{Thm}{Theorem}
\newtheorem{Con}{Conjecture}

\newtheorem{Lem}{Lemma}
\newtheorem{Def}{Definition}
\newtheorem{cor}{Corollary}

\newcommand{\qed}{\hfill $\Box$}

\begin{document}
\date{}
\maketitle 
{\def\thefootnote{}
}
\begin{abstract}
\noindent Let $a_n(k)$ be the $k$th coefficient of the $n$th cyclotomic polynomial $\Phi_n(x)$.
As $n$ ranges over the integers, $a_n(k)$ assumes only finitely many values.
For any such value $v$ we determine the density of integers $n$ such 
that $a_n(k)=v$. Also we study the average of the $a_n(k)$.
We derive analogous results for the $k$th Taylor coefficient of $1/\Phi_n(x)$ (taken
around $x=0$). We formulate various open problems.
\end{abstract}
\section{Introduction}
Let 
\begin{equation}
\label{definitie}
\Phi_n(x)=\prod_{j=1\atop (j,n)=1}^n(x-e^{2\pi ji\over n})=
\sum_{k=0}^{\varphi(n)}a_n(k)x^k,
\end{equation}
denote the $n$th cyclotomic polynomial and 
$\varphi$ Euler's totient function. 
If $k>\varphi(n)$, we put $a_n(k)=0$.
The coefficients $a_n(k)$ are integers.
In this paper we also consider the behaviour of the coefficients $c_n(k)$ in the
Taylor series of $1/\Phi_n(x)$ around $x=0$:
$${1\over \Phi_n(x)}=\sum_{k=0}^{\infty}c_n(k)x^k.$$
In the 19th century mathematicians were
already intrigued by the behaviour of $a_n(k)$, since the $a_n(k)$ seem to be amazingly
small. In a nutshell the history of the study of $a_n(k)$ can be described as inspired
by conjectures
about the $a_n(k)$ being small that were being proved wrong in due course. A 19th
century example being the conjecture that $a_n(k)\in \{-1,0,1\}$, which turned out to
be wrong once it was shown that $a_{105}(7)=-2$. A 21th century example being the
recent disproof of Gallot and Moree \cite{GM}
of the Beiter conjecture (dating back to 1968). This conjecture asserts that if $p<q<r$ 
are primes, then $|a_{pqr}(k)|\le (p+1)/2$. In \cite{GM} it is shown to be false for
every prime $p\ge 11$ and, moreover, that 
given any $\delta>0$ there exist infinitely many triples $(p_j,q_j,r_j)$
with $p_1<p_2<\ldots $ consecutive primes such that $|a_{p_jq_jr_j}(n_j)|>(2/3-\delta)p_j$
for $j\ge 1$.\\
\indent On noting that $x^n-1=\prod_{d|n}\Phi_d(x)$, we see that
$${x^n-1\over \Phi_n(x)}=\prod_{d|n,~d<n}\Phi_d(x)$$ 
is a polynomial of degree $n-\varphi(n)$ having
integer coefficients. {}From this we infer that the $c_n(k)$ are integers that only depend
on the congruence class of $k$ modulo $n$. Numerics suggest that 
the $c_n(k)$, like the $a_n(k)$, are surprisingly small. In constrast
to the $a_n(k)$ they only seem to have been studied as numbers of independent interest
in a paper by the second author \cite{Mor}.\\
\indent {}From the equality $x^n-1=\prod_{d|n}\Phi_d(x)$ one finds by inclusion and exclusion
that $\Phi_n(x)=\prod_{d|n}(x^d-1)^{\mu(n/d)}$, where $\mu$ denotes the M\"obius function.
On using that $\sum_{d|n}\mu(d)=0$ if $n>1$, we obtain, for $n>1$, 
\begin{equation}
\label{thanga}
\Phi_n(x)=\prod_{d|n}(1-x^d)^{\mu({n\over d})}.
\end{equation}
Sometimes it will be convenient to write
$$
\Phi_n(x)=\prod_{d=1}^{\infty}(1-x^d)^{\mu({n\over d})},
$$
where we have put $\mu(r)=0$ in case $r$ is not an integer.
Thus $a_n(k)$ is the coefficient of $x^k$ in $\prod_{d<k+1}(1-x^d)^{\mu({n\over d})}$.
Since $\mu\in \{-1,0,1\}$, we infer that $a_n(k)$ (and likewise $c_n(k)$) assumes only
finitely many values as $n$ ranges over the natural numbers.\\
\indent Put ${\cal A}(k)=\{a_n(k)|n\ge 1\}$
and ${\cal C}(k)=\{c_n(k)|n\ge 1\}$. 
Let $\chi_v(m)=1$ if $m=v$ and 0 otherwise.
Our main interest in this paper is to study these
sets and some associated quantities such as $$A(k):=\max\{|a_n(k)|:n\ge 1\}
{\rm ~and~}C(k):=\max\{|c_n(k)|:n\ge 1\}.$$ Also we are interested in the averages
$$M(a_n(k))=\lim_{x\rightarrow \infty}{\sum_{n\le x}a_n(k)\over x}{\rm ~and~}
M(c_n(k))=\lim_{x\rightarrow \infty}{\sum_{n\le x}c_n(k)\over x}.$$
Furthermore, for a given $v$ we consider the densities
$$\delta(a_n(k)=v)=\lim_{x\rightarrow \infty}{1\over x}\sum_{n\le x}\chi_v(a_n(k)){\rm ~and~}
\delta(c_n(k)=v)=\lim_{x\rightarrow \infty}{1\over x}\sum_{n\le x}\chi_v(c_n(k)).$$
We note that various authors have considered $A(k)$. In contrast we know of only one paper
dealing with $M(a_n(k))$ and sketching how $\delta(a_n(k)=v)$ can be computed. This is
a paper due to Herbert M\"oller \cite{HM}. In Section \ref{prev} we briefly discuss this
previous work and especially the latter paper as we will propose some improvements of it.\\
\hfil\break
\indent The results in this paper suggest that the behaviour of the $c_n(k)$ is so close
to that of the $a_n(k)$, making consideration of the $c_n(k)$ hardly worthwhile. This
is no longer true if we vary $k$ and keep $n$ fixed. E.g. one can have $|c_{pqr}(k)|=p-1$, 
whereas the estimate $|a_{pqr}(k)|\le 3p/4$ always holds (see \cite{Mor}). Furthermore,
in our average consideration the quantities $a_n(k)+c_n(k)$ and $a_n(k)-c_n(k)$ show up
in a natural way. In the M.Sc. thesis \cite{MH}, on which this paper is partly based, the
notion of $c_n(k)$ is not used, leading to slightly more complicated formulations of some of the
results.
\section{Basics and preliminaries}
Due to the fact that (\ref{thanga}) does not hold for $n=1$ various technical complications
arise. For this reason it turns out to be helpful to work with the following modified
cyclotomic coefficients (with $\epsilon =\pm 1$, by which
we mean $\epsilon\in \{-1,1\}$):
\begin{Def}
We let $a_n^{\epsilon}(k)$ be the $k$th Taylor coefficient (around $x=0$) in the product
$\prod_{d|n}(1-x^d)^{\epsilon\mu(n/d)}$, i.e., we have, for $|x|<1$,
\begin{equation}
\label{thangoo}
\prod_{d|n}(1-x^d)^{\epsilon\mu({n\over d})}=\sum_{k=0}^{\infty}a_n^{\epsilon}(k)x^k.
\end{equation}
\end{Def}
Sometimes we will use the latter identity in the following form (valid for $|x|<1$):
\begin{equation}
\label{thangoo1}
\prod_{d=1}^{\infty}(1-x^d)^{\epsilon\mu({n\over d})}=
\sum_{k=0}^{\infty}a_n^{\epsilon}(k)x^k.
\end{equation}
Note that the left hand side in (\ref{thangoo1}) equals $\Phi_n^{\epsilon}(x)$ in case $n>1$
and $-\Phi_1^{\epsilon}(x)$ in case $n=1$. 
{}From this we see that
$$a_n^1(k)=\cases{a_n(k) & if $n>1$;\cr -a_1(k) & if $n=1$,}{\rm ~and~}
a_n^{-1}(k)=\cases{c_n(k) & if $n>1$;\cr -c_1(k) & if $n=1$.}$$
\indent The basic properties of $\Phi_n(x)$ and its coefficients given below are quite useful. 
\begin{Lem}
\label{stard}
Let $q_1$ and $q_2$ be primes with $k<q_1<q_2$
and $(q_1q_2,n)=1$. Then $a_{nq_1}^{\epsilon}(k)=a_n^{-\epsilon}(k)$ 
and $a_{nq_1q_2}^{\epsilon}(k)=a_n^{\epsilon}(k)$.
\end{Lem}
{\it Proof}. An easy consequence of (\ref{thangoo}) and the properties of 
the M\"obius function.\qed\\

By $\gamma(n)=\prod_{p|n}p$ 
we denote the squarefree kernel of $n$.
\begin{Lem}
\label{tralala}
\hfil\break
{\rm 1)} We have $\Phi_{n}(x)=\Phi_{\gamma(n)}(x^{n/\gamma(n)})$.\\
{\rm 2)} We have $\Phi_{2n}(x)=\Phi_n(-x)$ if $n>1$ is odd.\\ 
{\rm 3)} We have $x^{\varphi(n)}\Phi_n(1/x)=\Phi_n(x)$ if $n>1$.
\end{Lem}
\indent In terms of the coefficients
the three properties of Lemma \ref{tralala} imply (respectively):
\begin{equation}
\label{naarkwadraatvrijekern}
a_n^{\epsilon}(k)=\cases{a_{\gamma(n)}^{\epsilon}({k\gamma(n)\over n}) &if ${n\over \gamma(n)}|k$;\cr
0 & otherwise,}
\end{equation}
\begin{equation}
\label{verdubbeling}
a_{2n}^{\epsilon}(k)=(-1)^ka_n^{\epsilon}(k)~{\rm ~if~}2\nmid n;
\end{equation}
\begin{equation}
\label{symm}
a_{n}^{\epsilon}(k)=\epsilon a_n^{\epsilon}(\varphi(n)-k)~{\rm ~for~}n>1,~0\le k\le \varphi(n),
\end{equation}
where to prove (\ref{verdubbeling}) in case $n=1$ we used that $\Phi_2^{\epsilon}(x)=-\Phi_1^{\epsilon}(-x)$.
In order to prove (\ref{symm}) in case $\epsilon=-1$, we used the additional observation that
$${1-x^n\over \Phi_n(x)}=\sum_{k=0}^{n-\varphi(n)}a_n^{-1}(k)x^k.$$
It is not difficult to derive Lemma \ref{tralala} from (\ref{thanga}),
see e.g. Thangadurai \cite{Thanga}.\\
\indent Note that, for $|x|<1$, we have
\begin{equation}
\label{lehmertje}
\prod_{d=1}^{\infty}(1-x^d)^{\mu({n\over d})}=\prod_{d=1}^{\infty}\left(1-\mu({n\over d})x^d+
{1\over 2}\mu({n\over d})(\mu({n\over d})-1)\sum_{j=2}^{\infty}x^{jd}\right),
\end{equation}
where we used the observation that, for $|x|<1$,
\begin{equation}
\label{thanga2}
(1-x^d)^{\mu({n\over d})}=1-\mu({n\over d})x^d+
{1\over 2}\mu({n\over d})(\mu({n\over d})-1)\sum_{j=2}^{\infty}x^{jd}.
\end{equation}
{}From (\ref{lehmertje}) it is not difficult to derive a formula for $a^{1}_n(k)$ for a fixed $k$;
this is just the coefficient of $x^k$ in the right hand side of (\ref{lehmertje}) (this approach
seems to be due to D.H. Lehmer \cite{DL}). We thus obtain, 
$$\cases{
a^{1}_n(1)=-\mu(n);\cr
a^{1}_n(2)=\mu(n)^2/2-\mu(n)/2-\mu(n/2);\cr
a^{1}_n(3)=\mu(n)^2/2-\mu(n)/2+\mu(n/2)\mu(n)-\mu(n/3).}$$
More generally, we have
$$
a_n^{1}(k)=\sum c(k_1,...,k_s;e_1,\dots,e_s)\mu({n\over k_1})^{e_1}\cdots \mu({n\over k_s})^{e_s},
$$
where the sum is over all partitions $k_1+\ldots+k_s$ of all the
integers $\le k$ with $k_1\ge k_2\ge \cdots \ge k_s$ and over all $e_1,\dots,e_s$ with
$1\le e_j\le 2$ for $1\le j\le s$. The terms in (\ref{uitdrukking}) for which
$e_1+\ldots+e_s$ is even we add together to obtain $\alpha_n(k)$, the {\it even part} 
of
$a_n^{1}(k)$. Similarly, we group the terms with $e_1+\ldots+e_s$ odd together, to form
the {\it odd part}, $\beta_n(k)$,  of $a_n^{1}(k)$. (To the authors knowledge the 
even and odd part of $a_n^1(k)$ have not
been defined and considered before.) 
For example, $\alpha_n(2)=\mu(n)^2/2$ and $\beta_n(2)=-\mu(n)/2-\mu(n/2)$. Note that
\begin{eqnarray}
\label{uitdrukking}
a_n^{\epsilon}(k)&=&\sum c(k_1,...,k_s;e_1,\dots,e_s)(\epsilon \mu ({n\over k_1}))^{e_1}\cdots 
(\epsilon \mu({n\over k_s}))^{e_s}.\cr
&=&\alpha_n(k) + \epsilon \beta_n(k).
\end{eqnarray}
We have
$\alpha_n(k)=(a_n^1(k)+a_{n}^{-1}(k))/2$ and $\beta_n(k)=(a_n^1(k)-a_{n}^{-1}(k))/2$. In particular, 
$2\alpha_n(k),2\beta_n(k)\in \mathbb Z$. 
{}From (\ref{uitdrukking}) and the properties of the M\"obius
function it follows that if $p$ and $q$ are two distinct primes exceeding $k$ 
with $(pq,n)=1$, then $\alpha_{pn}(k)=\alpha_n(k)$, $\beta_{pn}(k)=-\beta_n(k)$, 
$\alpha_{pqn}(k)=\alpha_n(k)$ and $\beta_{pqn}(k)=\beta_n(k)$. The reason, as we will see, for
distinguishing between the odd and even part, is that the odd part does not contribute to
te average, i.e. $M(\beta_n(k))=0$.\\
\indent The {\it Ramanujan sum} $r_n(m)$ is defined by
\begin{equation}
\label{ramadefinitie}
r_n(m)=\sum_{1\le k\le n\atop (k,n)=1}e^{2\pi i mk\over n}=\sum_{1\le k\le n\atop (k,n)=1}\zeta_n^{mk}.
\end{equation}
Alternatively one can write $r_n(m)={\rm Tr}_n(\zeta_n^m)$, where by Tr$_n$ we denote the trace over
the cyclotomic field $\mathbb Q(\zeta_n)$. It follows at once from the properties of the trace
that $r_n(m)=r_n((n,m))$. Since $\zeta_n^m$ is an algebraic integer, it follows that 
$r_n(m)$ is an integer.\\ 
\indent The Ramanujan sums have many properties of which we will need only the following two.
\begin{Lem}
\label{basicramanujan}
We have
$r_n(m)=\sum_{d|(n,m)}d\mu({n\over d})$ and
$$r_n(m)=\mu\left({n\over (n,m)}\right){\varphi(n)\over \varphi({n\over (n,m)})}.$$
\end{Lem}
Nicol \cite{nicol} showed that Ramanujan sums and cyclotomic polynomials are closely related, by establishing
that
$$\Phi_n(x)=\exp\Big(-\sum_{m=1}^{\infty}{r_n(m)\over m}x^m\Big)~{\rm and~}\sum_{m=1}^nr_n(m)x^{m-1}=
(x^n-1){\Phi_n'(x)\over \Phi_n(x)}.$$
\subsection{Some sums involving the M\"obius function}
In order to evaluate $M(a_n(k))$ we will need to evaluate $\sum_{m\le x,~(m,r)=1}\mu(m)^k$
with $1\le k\le 2$. That is done in Lemmas \ref{mobiuseen} and \ref{mobiustwee}.
\begin{Lem}
\label{mobiuseen}
Let $r\ge 1$ be an integer. We have
$$\sum_{m\le x\atop (m,r)=1}\mu(m)^2={6x\over \pi^2\prod_{p|r}(1+{1\over p})}+O(\sqrt{x}\varphi(r)),$$
where the implied constant is absolute.
\end{Lem}
{\it Proof}. We have, by inclusion and exclusion,
$$\sum_{m\le x\atop (m,r)=1}\mu(m)^2=\sum_{d\le \sqrt{x}\atop (d,r)=1}\mu(d)A_r({x\over d^2}),$$
where $A_r(x)$ denotes the number of integers $n\le x$ that are coprime with $r$. Note that
$$\Big[{x\over r}\Big]\varphi(r)\le A_r(x)\le \Big[{x\over r}\Big]\varphi(r)+\varphi(r)$$
and hence $A_r(x)=\varphi(r)x/r+O(\varphi(r))$. On using the latter estimate we obtain
\begin{eqnarray}
\sum_{m\le x\atop (m,r)=1}\mu(m)^2&=&x{\varphi(r)\over r}\sum_{d\le \sqrt{x}\atop (d,r)=1}{\mu(d)\over d^2}
+O(\sqrt{x}\varphi(r)).\nonumber\\
&=&x{\varphi(r)\over r}\sum_{(d,r)=1}^{\infty}{\mu(d)\over d^2}+O(\sqrt{x}\varphi(r)).\nonumber\\
&=&{6x\over \pi^2\prod_{p|r}(1+{1\over p})}+O(\sqrt{x}\varphi(r)),\nonumber
\end{eqnarray}
where we used that
$${\varphi(r)\over r}\sum_{(d,r)=1}^{\infty}{\mu(d)\over d^2}={\varphi(r)\over r}
\prod_{p\nmid r}(1-{1\over p^2})={\varphi(r)\over \zeta(2)r\prod_{p|r}(1-{1\over  p^2})}
={1\over \zeta(2)\prod_{p|r}(1+{1\over p})}$$
and $\zeta(2)=\pi^2/6$. \qed
\begin{Lem}
\label{mobiustwee}
We have $$\sum_{m\le x\atop (m,r)=1}\mu(m)=o(x),$$ where the implied constant may depend on $r$.
\end{Lem}
To prove the lemma, we will
apply the Wiener-Ikehara Tauberian theorem in the following form.
\begin{Thm}
\label{vertaubung}
Let $f(s)=\sum_{n=1}^{\infty}a_n/n^s$ be a Dirichlet series. Suppose there exists a Dirichlet series
$F(s)=\sum_{n=1}^{\infty}b_n/n^s$ with positive real coefficients such that\\
{\rm (a)} $|a_n|\le b_n$ for all $n$;\\
{\rm (b)} the series $F(s)$ converges for Re$(s)>1$;\\
{\rm (c)} the function $F(s)$ can be extended to a meromorphic function in the region Re$(s)\ge 1$ having
no poles except for a simple pole at $s=1$.\\
{\rm (d)} the function $f(s)$ can be extended to a meromorphic function in the region Re$(s)\ge 1$ having
no poles except possibly for a simple pole at $s=1$ with residue $r$.\\
Then
$$\sum_{n\le x}a_n=rx+o(x),~x\rightarrow \infty.$$
In particular, if $f(s)$ is holomorphic at $s=1$, then $r=0$ and
$\sum_{n\le x}a_n=o(x)$ as $x\rightarrow \infty$.
\end{Thm} 
{\it Proof of Lemma} \ref{mobiustwee}. We apply the Wiener-Ikehara theorem with $F(s)=\zeta(s)$ and
$$f(s)=\sum_{(n,r)=1}{\mu(n)\over n^s}={1\over \zeta(s)\prod_{p|r}(1-{1\over p^s})}.$$
Of course $F(s)$ satisfies the required properties and has a simple pole at $s=1$ with residue one.
Since the finite product in the formula for $f(s)$ is regular for Re$(s)>0$, the result follows
on using the well-known fact that $1/\zeta(s)$ can be extended to a meromorphic function in
the region Re$(s)\ge 1$ (and hence $r=0$). This completes the proof.\qed\\

\noindent Landau \cite[\S 173 \& \S 174]{L} gave estimates for
$\sum_{n\le x,~n\equiv l({\rm mod~}k)}\mu(n)^r$, with $1\le r\le 2$. Using these, an
alternative proof of Lemmas \ref{mobiuseen} and \ref{mobiustwee} can be given.\\
\indent The following result is found on combining Lemma \ref{mobiuseen} and \ref{mobiustwee}.
\begin{Lem}
\label{mobiusdrie}
Let $\epsilon=\pm 1$ and $r\ge 1$. We have, as $x$ tends to infinity,
$$\sum_{m\le x,~\mu(m)=\epsilon\atop (m,r)=1}1\sim {3x\over \pi^2\prod_{p|r}(1+{1\over p})}.$$
\end{Lem}

\section{Previous work}
\label{prev}
In this section we discuss previous work on the quantities defined in the introduction. Of
those $A(k)$ has received quite a bit of attention. H. M\"oller \cite{HM} gave a table
for $A(k)$ for $1\le k\le 20$ which we reproduce below (this before Endo \cite{Endo}, 
who proved that $A(k)=1$ for $k\le 6$).\\

\centerline{{\bf Table 1:} $A(k)$ for $1\le k\le 30$}
\begin{center}
\begin{tabular}{|c|c|c|c|c|c|c|c|c|c|c|c|c|c|c|c|c|}
\hline
$k$  & $1$ & $2$ & $3$ & $4$ & $5$ & $6$ & $7$ & $8$ & $9$ & $10$ & $11$ & $12$ & $13$ & $14$ & $15$ \\
\hline
$A(k)$  & $1$ & ${1}$ & ${1}$ & ${1}$ & ${1}$ & 
$1$ & $2$ & $1$ & ${1}$ & $1$ & $2$ & $1$ & $2$ & $2$ & $2$\\
\hline
$k$ & $16$ & $17$ & $18$ & $19$ & $20$ & $21$ & $22$ & $23$ & $24$ & $25$ & $26$ & $27$ & $28$ & $29$ & $30$\\
\hline
$A(k)$ & $2$ & $3$ & $3$ & $3$ & $3$ & $3$ &  $3$ & $4$ & $3$ & $3$ & $3$ & $3$ & $4$ & $4$ & $5$ \\
\hline
\end{tabular}
\end{center}
\medskip
The table suggests that $A(k)\le k$ for every $k\ge 1$, this is however very far from
being the case as given $r>1$ we have $A(k)>k^r$ for all $k$ sufficiently
large as M\"oller proved.
Bachman \cite{gennady} extended work of several earlier authors (see
the references he gives) and established the best result to
date stating that
\begin{equation}
\label{gennady}
\log A(k)=C_0{\sqrt{k}\over (\log k)^{1/4}}\left(1+O\left({\log \log k\over \sqrt{\log k}}\right)\right),
\end{equation}
where $C_0$ can be explicitly given.\\
\indent In M\"oller's approach $a_n(k)$ is 
connected with partitions of $k$ through the basic
formula (\ref{simpelzegA}). A partition can be
identified with a sequence $\{n_j\}_{j=1}^{\infty}$ of non-negative integers satisfying
$\sum_{j}jn_j=m$. 
Without loss of generalisation we can denote a partition, $\lambda$, of $k$ as
$\lambda=(k_1^{n_{k_1}}\cdots k_s^{n_{k_s}})$, where
$n_{k_1}\ge n_{k_2}\ge \dots \ge n_{k_s}\ge 1$ (thus the number $k_j$ occurs $n_{k_j}$ times
in the partition). The set of all partitions of $m$ will be denoted by ${\cal P}(m)$. 
The number of different partitions of $m$ is denoted by $p(m)$. Hardy and
Littlewood in 1918, and Uspensky independently in 1920, proved that
\begin{equation}
\label{harlit}
p(m)\sim {e^{\pi \sqrt{2m/3}}\over 4m\sqrt{3}}~{\rm ~as~}m\rightarrow \infty.
\end{equation}
In case $\epsilon=1$ and $n>1$ the following result is \cite[Satz 2]{HM}.
\begin{Lem} 
\label{simpelzeg}
For $n\ge 1$, $k\ge 0$ and $\epsilon=\pm 1$, we have
\begin{equation}
\label{simpelzegA}
a_n^{\epsilon}(k)=\sum_{\lambda=(k_1^{n_{k_1}}\dots k_s^{n_{k_s}})\in {\cal P}(k)\atop
n_{k_1}\ge \dots \ge n_{k_s}\ge 1}\prod_{j=1}^s(-1)^{n_{k_j}}\left({\epsilon\mu({n\over k_j})\atop n_{k_j}}\right),
\end{equation}
where the sum is over all partitions $\lambda$ of ${\cal P}(k)$.
\end{Lem}
{\it Proof}. The Taylor series of $(1-x)^a$ equals, for $|x|<1$, 
$\sum_{j=0}^{\infty}(-1)^j({a\atop j})x^j$, where
$({a\atop j})=a(a-1)\cdots (a-(j-1))/j!$. Using this we infer that
\begin{equation}
\label{nogeentje}
(1-x^d)^{\epsilon\mu({n\over d})}=\sum_{j=0}^{\infty}(-1)^j\left({\epsilon\mu({n\over d})\atop j}\right)x^{dj},~|x|<1,
\end{equation}
The proof now follows from (\ref{nogeentje}) and (\ref{thangoo1}). \qed\\

\noindent Earlier D. Lehmer \cite{DL} had used a formula for $a_n(k)$, expressing it
in terms of Ramanujan sums. The formula above turns out to be more practical.
Our proof above shows that formula (\ref{simpelzegA})  is a triviality, whereas M\"oller's 
ingenious and rather involved proof of it (he considered only $\epsilon=1$) obscures this.\\
\indent Comparison of (\ref{thanga2}) and (\ref{nogeentje}) yields
\begin{equation}
\label{veertien}
(-1)^j\left({\epsilon\mu({n\over d})\atop j}\right)=\cases{1 & if $j=0$;\cr
-\epsilon\mu(n/d) &if $j=1$;\cr
\epsilon\mu(n/d)(\epsilon\mu(n/d)-1)/2 &if $j\ge 2$.}
\end{equation}
Thus the product appearing in (\ref{simpelzegA}) is in $\{-1,0,1\}$ and
it follows from (\ref{simpelzegA}) that $|a_n^{\epsilon}(k)|\le p(k)$. By the asymptotic formula for $p(m)$ given above it then
follows that $\log A(k)\ll \sqrt{k}$ is the trivial upper bound for $\log A(k)$.
{}From Lemma \ref{simpelzeg} M\"oller infers that even $|a_n(k)|\le p(k)-p(k-2)$.
To see this note that the partitions having 1 occurring at least twice do not contribute if
$\mu(n)\in \{0,1\}$. If $\mu(n)=-1$, then either $\mu(2n)=1$ or $\mu(n/2)=1$. This in
combination with (\ref{verdubbeling}) allows us then to argue as before and leads
us to the same bound. Similarly we have $|c_n(k)|\le p(k)-p(k-2)$.\\ 
\indent M\"oller uses Lemma \ref{simpelzeg} to show that
$$M(a_n(k))=\lim_{x\rightarrow \infty}{\sum_{n\le x}a_n(k)\over x}$$
exists and gives a formula for it. To do so
he first 
computes the Dirichlet series $D_k(s):=\sum_{n=1}^{\infty}a_n(k)n^{-s}$. 
The required average is then $\lim_{s\rightarrow \infty}(s-1)D_k(s)$. The
expression so obtained is rather complicated and requires work to be simplified.
Here we rederive his result (see Lemma \ref{halfweg}) in a more direct way by simply evaluating the
averages 
$$\lim_{x\rightarrow \infty}{1\over x}\sum_{n\le x}\prod_{j=1}^s (-1)^{n_{k_j}}\left({\mu({n\over k_j})\atop
n_{k_j}}\right)$$
and then summing over the partitions of $k$.\\
\indent M\"oller's result shows that  
$M(a_n(k))=6e_k/\pi^2$, with $e_k$  a rational 
number.  For $1\le k\le 20$ we give the value of $e_k$ in Table 2 (our table agrees with
the one given in \cite{HM}, except for the incorrect values $e_{10}=319/1440$ 
and $e_{16}=733/2016$ appearing there). \\

\newcommand\T{\rule{0pt}{2.6ex}}
\newcommand\B{\rule[-1.4ex]{0pt}{0pt}} 
\centerline{{\bf Table 2:} Scaled average, $e_k=\zeta(2)M(a_n(k))$, of $a_n(k)$ }
\begin{center}
\begin{tabular}{|c|c|c|c|c|c|c|c|c|c|c|c|}
\hline
\T \B $k$  & $1$ & $2$ & $3$ & $4$ & $5$ & $6$ & $7$ & $8$ & $9$ & $10$ \\
\hline
\T \B $e_k$  & $0$ & ${1\over 2}$ & ${1\over 6}$ & ${1\over 3}$ & ${1\over 8}$ & 
${7\over 24}$ & ${1\over 18}$ & ${7\over 24}$ & ${19\over 144}$ & ${31\over 160}$\\
\hline
\T \B $k$  & $11$ & $12$ & $13$ & $14$ & $15$ & $16$ & $17$ & $18$ & $19$ & $20$\\
\hline
\T \B $e_k$  & ${1\over 16}$ & ${55\over 192}$ & ${13\over 288}$ & ${61\over 288}$ & ${2287\over 20160}$ &
${733\over 4032}$ & ${667\over 8064}$ & ${79\over 336}$ & ${55\over 1344}$ & ${221\over 960}$  \\
\hline
\end{tabular}
\end{center}
\medskip
Regarding $e_k$
M\"oller proposed:
\begin{Con} {\rm \cite{HM}}. 
\label{vermoedentwee}
Let $k\ge 1$. Write $M(a_n(k))=6e_k/\pi^2$.\\
{\rm 1)} We have $0\le e_k\le 1/2$.\\
{\rm 2)} We have $(-1)^k(e_k-e_{k+1})>0$ ("see-saw conjecture").
\end{Con}
M\"oller stated that with help of an IBM 7090 he wanted
to check his conjecture for further values of $k$. Had he carried this out, he would have
discovered that $(-1)^{34}(e_{34}-e_{35})=-18059/4626720<0$. Other
counterexamples occur at $k=35,45$ and $94$. Indeed, we would not be surprised if part 2 of
the Conjecture is violated for infinitely many $k$. The see-saw conjecture, if true, would
imply that $\sum_{k=1}^m (-1)^k(e_k-e_{k+1})>0$ for every $m\ge 1$. The truth of the latter assertion
is still open.\\
\indent On the other hand, part 1 of the Conjecture is true for $k\le 100$. The numbers
$e_k$ seem to be decreasing to zero and their size seems to be related to the number
of prime factors of $k$, the more prime factors the larger $e_k$ seems to be.\\
\indent As already pointed out by M\"oller one could use his method to study the value distribution of
$a_n(k)$ in case e.g. $A(k)=1$ by considering the integer $a_n(k)(a_n(k)-1)/2$ to determine
$\delta(a_n(k)=-1)$ for example. This then yields a sum with $p(k)^2+p(k)$ terms and this results
in an algorithm that has worse complexity than that provided by Theorem \ref{vier} below. Aside from this,
this seems to be, from the practical point of view, an unwieldy method. A more
practical method will be presented in Section \ref{vvvier}.\\
\indent A further result which is of relevance to us, is the following one.
\begin{Thm}
\label{z}
Let $m\ge 1$ be an integer and $\epsilon=\pm 1$.
Then $$\{a_{mn}^{\epsilon}(k):n\ge 1,~k\ge 0\}=\mathbb Z.$$
\end{Thm}
For a proof and the prehistory of this result see Ji, Li and Moree \cite{JLM}.

\section{Computation of ${\cal A}(k)$ and ${\cal C}(k)$}
Recall that ${\cal A}(k)=\{a_n(k)|n\ge 1\}$
and ${\cal C}(k)=\{c_n(k)|n\ge 1\}$. Throughout this section we assume that $k\ge 1$.
\begin{Lem} 
\label{brovo}
We have 
$$\{-1,0,1\}\subseteq \{a_n(k):n>1\}{\rm ~and~}\{-1,0,1\}\subseteq \{c_n(k):n>1\}.$$
\end{Lem}
{\it Proof}. In formula (\ref{simpelzeg}) there is always
the term $-\epsilon\mu(n/k)$. Let us take $n=ck\prod_{p\le k}p$, where $c$ only has
prime divisor $>k$. Then all the terms of the form $\mu(n/r)$
with $1\le r<k$ are zero (since either $r\nmid n$ or $n/r$ is not squarefree) and we obtain
that $a_n^{\epsilon}(k)=-\epsilon\mu(c)(-1)^{\pi(k)}$, where $\pi(x)$ as usual denotes the number of primes $p\le x$
not exceeding $x$. In particular, it follows that $a_n^{\epsilon}(k)$ always assumes the values $-1,0$ and $1$. 
Since $n>1$ for these examples, $a_n^{\epsilon}(k)=a_n(k)$ if $\epsilon=1$ and
equals $c_n(k)$ if $\epsilon=-1$, and the result follows. \qed
\begin{Lem} 
\label{schuif}
We have ${\cal A}(k)=\{a_n^1(k)|n\ge 1\}$ and ${\cal C}(k)=\{a_n^{-1}(k)|n\ge 1\}$.
\end{Lem}
{\it Proof}. Using Lemma \ref{brovo} one infers that 
$${\cal A}(k)=\{a_1(k)\}\cup \{a_n(k)|n>1\}
=\{a_n(k)|n>1\}=\{a_n^1(k)|n>1\}.$$ Likewise we find ${\cal C}(k)=\{a_n^{-1}(k)|n\ge 1\}$. \qed\\

The next lemma follows on applying the latter lemma in combination with Lemma \ref{stard}. 
\begin{Lem} 
\label{simcard}
We have ${\cal C}(k)={\cal A}(k)$.
\end{Lem}
Lemma \ref{naareindig} allows one to deduce that ${\cal A}(k)$ is a finite set.
\begin{Lem}
\label{naareindig}
Put $N_k={\rm lcm}(1,2,\cdots,k)\prod_{p\le k}p$. We can uniquely decompose $n$
as $n=n_kc_k$ with $(c_k,N_k)=1$
and $n_k$ and $c_k$ natural numbers. Let $\epsilon=\pm 1$.
There exist functions $A_1$ and $B_1$ with as domain the divisors of $N_k$ such that
$$a_n^{\epsilon}(k)=\cases{A_1(n_k)\mu(c_k)^2+\epsilon B_1(n_k)\mu(c_k) &if $n_k|N_k$;\cr
0 &otherwise.}$$
\end{Lem}
{\it Proof}. The assertion regarding the uniqueness of the decomposition $n=n_kc_k$ is trivial. 
For a given $n$ only those partitions $k_1,k_2,\ldots,k_s$ contribute to 
(\ref{uitdrukking}) for which
$n/k_i$ is an integer for $1\le i\le s$. Note that $k_i|n_k$. Thus, we can write
$$\mu({n\over k_1})^{e_1}\cdots \mu({n\over k_s})^{e_s}=
\mu({n_k\over k_1})^{e_1}\cdots \mu({n_k\over k_s})^{e_s}\mu(c_k)^{e_1+\dots+e_s}.$$
If $n_k\nmid N_k$, then none of the integers $n_k/k_1,...,n_k/k_s$ is squarefree and 
$a_n^{\epsilon}(k)=0$, so assume that $n_k|N_k$. 
On using that $\mu(r)^w$ with $w\ge 1$ either equals $\mu(r)$ or $\mu(r)^2$, it follows
from (\ref{uitdrukking}) that $a_n^1(k)=A_1(n_k)\mu(c_k)^2+B_1(n_k)\mu(c_k)$.\\ 
\indent Note that $\alpha_n(k)=A_1(n_k)\mu(c_k)^2$ and that
$\beta_n(k)=B_1(n_k)\mu(c_k)$. Thus $a_n^{-1}(k)=\alpha_n(k)-\beta_n(k)=
A_1(n_k)\mu(c_k)^2-B_1(n_k)\mu(c_k)$. \qed\\ 

\noindent Using formula (\ref{naarkwadraatvrijekern}) it is seen that in the latter lemma $N_k$ can
be replaced by $k\prod_{p\le k}p$.
\begin{Lem}
\label{naareindig2}
Lemma {\rm \ref{naareindig}} holds true also with $N_k$ replaced by $M_k=k\prod_{p\le k}p$.
\end{Lem}
\noindent The above lemma 
again shows that ${\cal A}(k)$ is a finite set, thus if
we fix $k$, there are only finitely many possibilities for the values of the
coefficient of $x^k$ in a cyclotomic polynomial.\\
\indent The latter lemma in combination with Lemma \ref{schuif} yields the following result. 
\begin{Lem} 
\label{valueset}
Let $M_k=k\prod_{p\le k}p$. Then
$\{0,a_{d}^1(k),a_d^{-1}(k)~|~d|M_k\}={\cal A}(k)$.
\end{Lem}
We have $|a_n(k)|\le \max_{n\ge 1}|a_n(k)|=A(k)$. See Table 1 for the values
of $A(k)$ for $1\le k\le 30$. Define ${\cal A}^{\epsilon}(k)=\{a_d^{\epsilon}(k)~|~d|M_k\}$.
Numerics show that mostly
${\cal A}^{\epsilon}(k)={\cal A}(k)$. If ${\cal A}^{\epsilon}(k)$ is strictly
included in ${\cal A}(k)$ (which happens for example for $k=48,54$, $\epsilon=\pm 1$ and 66, $\epsilon=-1$), then for 
the $k$ for which we did the computation ($k\le 73$),
the set $\{a_d^{\epsilon}(k)~|~d|M_k\}$ equals ${\cal A}(k)$ with one element omitted and this element is either
$A(k)-1$ or $-A(k)+1$.\\
\indent We next show that the inclusion in Lemma \ref{brovo} is strict for $k\ge 13$. Our
proof rests on 
the following rather elementary result on prime numbers. 
\begin{Lem}
\label{bertie}
For $k\ge 13$ there are consecutive odd primes $p_1<p_2<p_3$ such that $p_3\le k<p_1+p_2$.
\end{Lem}
{\it Proof}. Breusch \cite{Breusch} proved that for $x\ge 48$ there is at least one prime in $[x,9x/8]$ 
(this strengthens {\it Bertrand's Postulate} 
asserting that there is always a prime between $x$ and
$2x$, provided $x\ge 2$).
Let $\alpha=1.32$. A little computation shows that the above result implies that for $x\ge 9$ there
is at least one prime in $[x,\alpha x]$. One checks that the assertion is true for $k\in [13,21)$.
Assume that $k\ge 21~(\ge 9\alpha^3)$. Let $p_3$ be the largest prime not exceeding $k$ and let
$p_1$ and $p_2$ be primes such that $p_1,p_2$ and $p_3$ are consecutive primes. Then
$p_3\ge k/\alpha$, $p_2\ge k/\alpha^2$ and $p_1\ge k/\alpha^3$. On noting that
$p_1+p_2\ge k(1/\alpha+1/\alpha^2)>k$, the proof is then completed. \qed

\begin{Lem}
For $k\ge 13$ we have $\{-2,-1,0,1\}\in {\cal A}(k)$ (and thus $A(k)\ge 2$).
\end{Lem}
{\it Proof}. Let $p_1,p_2$ and $p_3$ be odd primes satisfying the condition of Lemma \ref{bertie}. 
Using  (\ref{thanga}) we infer that
$$\Phi_{p_1p_2p_3}(x)\equiv {(1-x^{p_3})\over (1-x)}(1-x^{p_1})(1-x^{p_2})\equiv (1+x+\cdots+x^{p_3-1})
(1-x^{p_1}-x^{p_2}),$$
where we computed modulo $x^{k+1}$.
It follows that $a_{p_1p_2p_3}(k)=-2$.\qed

The following result shows that if $k$ is odd, then ${\cal A}(k)$ is symmetric, that is
if $v\in {\cal A}(k)$, then also $-v\in  {\cal A}(k)$.
\begin{Lem}
\label{symmie}
If $k$ is odd, then ${\cal A}(k)=-{\cal A}(k)$, that is ${\cal A}(k)$ is symmetric.
\end{Lem}
{\it Proof}. Assume that $v\in {\cal A}(k)$. If $v=0$ there is nothing to prove, so assume that $v\ne 0$. Since $M_k$ is odd,
it follows by Lemma \ref{valueset} that $a_d^{\epsilon}(k)=v$ for some odd integer $d$ and
$\epsilon \in \{-1,1\}$. Then, by (\ref{verdubbeling}), 
$a_{2d}^{\epsilon}(k)=(-1)^ka_d^{\epsilon}(k)=-v$. On invoking Lemma \ref{simcard}, the proof
is then completed.\qed

\subsection{Numerical evaluation of $a_n(k)$ and $c_n(k)$ for small $k$ }
For our purposes it is relevant to be able to 
numerically evaluate $a_n(k)$ for small $k$ and large $n$. A computer package like Maple
evaluates $a_n(k)$ by evaluating the whole polynomial $\Phi_n(x)$. For large $n$ this
is far too costly. Instead it is more efficient to use (\ref{thanga}) and expand for
every $d$ for which $\mu(n/d)\ne 0$, $(1-x^d)^{\mu(n/d)}$ as a Taylor series up to $O(x^{k+1})$
and multiply all these series together.\\
\indent  The most efficient method to date to compute $a_n(k)$ for small $k$ is due to
Grytczuk and Tropak \cite{GT}. First they apply formula (\ref{naarkwadraatvrijekern}). Thus
it is enough to compute $a_n(k)$ with $n$ squarefree. If $\phi(n)<k$, then $a_n(k)=0$, so
we may assume that $\phi(n)\ge k$. Let $d=(n,\prod_{p\le k}p)$. Put
$T_r=\mu(n)\mu((r,d))\varphi((r,d))$. Compute $b_0,\dots,b_k$ recursively by $b_0=1$ and
$$b_j=-{1\over j}\sum_{m=0}^{j-1}b_mT_{j-m}{\rm ~for~}1\le j\le k.$$
Then $b_k=a_n(k)$. Their proof uses the formula
\begin{equation}
\label{tropakje}
a_n(k)=-{1\over k}\sum_{m=0}^{k-1}a_n(m)r_n(k-m)~{\rm for~}k\ge 1,
\end{equation}
which follows by Vi\`ete's and Newton's formulae from (\ref{definitie}) and it uses the
second formula
of Lemma \ref{basicramanujan}. However, an alternative proof of (\ref{tropakje}) is obtained
on using the following observation together with the first formula
of Lemma \ref{basicramanujan}.
\begin{Lem}
\label{achtt}
Suppose that, as formal power series, 
$$\prod_{d=1}^{\infty}(1-x^d)^{-a_d}=\sum_{d=0}^{\infty}r(d)x^d,$$
then $dr(d)=\sum_{j=1}^d r(d-j)\sum_{k|j}ka_k$.
\end{Lem}
{\it Proof}. Taking the logarithmic derivative of $\prod_{d=1}^{\infty}(1-x^d)^{-a_d}$ we obtain
$${\sum_{d=1}^{\infty}dr(d)x^d\over \sum_{d=0}^{\infty}r(d)x^d}=x{d\over dx}\log 
\prod_{d=1}^{\infty}(1-x^d)^{-a_d}=\sum_{j=1}^{\infty}(\sum_{k|j}ka_k)x^j,$$
whence the result follows. \qed\\

\noindent The following result  generalizes the Grytczuk-Tropak algorithm to the efficient
computation of $a_n^{\epsilon}(k)$ for small $k$ and large $n$.
\begin{Lem}
\label{geetee}
Let $n$ be squarefree and put $d=(n,\prod_{p\le k}p)$. Furthermore, we put
$T_r=\mu(n)\mu((r,d))\varphi((r,d))$. Compute $b_0,\dots,b_k$ recursively by $b_0=1$ and
$$b_j=-{\epsilon\over j}\sum_{m=0}^{j-1}b_mT_{j-m}{\rm ~for~}1\le j\le k.$$
Then $a_n^{\epsilon}(k)=b_k$.
\end{Lem}
{\it Proof}. For $k=0$ we have $1=b_0=a_n^{\epsilon}(k)$ and so we may assume that $k\ge 1$.
Apply Lemma \ref{achtt} with $a_d=-\epsilon \mu({n\over d})$ (thus $r(d)=a_n^{\epsilon}(d)$).
We obtain by part 1 of Lemma \ref{basicramanujan} that 
\begin{equation}
\label{tropakje2}
a_n^{\epsilon}(k)=-{\epsilon\over k}\sum_{m=0}^{k-1}a_n^{\epsilon}(m)r_n(k-m)~{\rm for~}k\ge 1,
\end{equation}
The proof now follows if we show that $r_n(r)=\mu(n)\mu((r,d))\varphi((r,d))$ for $1\le r\le k$.
Since by assumption $n$ is squarefree, part 2 of Lemma \ref{basicramanujan} implies
that $r_n(r)=\mu(n)\mu((n,r))\varphi((n,r))$. In case $1\le r\le k$ this can be rewritten as
\begin{eqnarray}
r_n(r)&=&\mu(n)\mu((n,r,\prod_{p\le k}p)\varphi((n,r,\prod_{p\le k}p))\nonumber\cr
&=& \mu(n)\mu((d,r))\varphi((d,r)).
\end{eqnarray}
Thus the proof is completed. \qed\\

\noindent {\tt Algortihm to compute $a_n^{\epsilon}(k)$}. If $n$ is not squarefree, then apply 
(\ref{naarkwadraatvrijekern}). Thus we may assume that $n$ is squarefree.\\
{\it The case} $\epsilon=1$. If $n>\varphi(k)$, then $a_n^1(k)=0$, otherwise compute $a_n^1(k)$ using
Lemma \ref{geetee}.\\
{\it The case} $\epsilon=-1$. Let $0\le k_1<n$ be such that $k_1\equiv k({\rm mod~}n)$. Then
$a_n^{-1}(k)=a_{n}^{-1}(k_1)$. If $k_1>n-\varphi(n)$, then $a_n^{-1}(k)=0$, otherwise
compute $a_n^{-1}(k)$ using Lemma \ref{geetee}.\\

\noindent In case $n>1$, $a_n^1(k)=a_n(k)$ and the above algorithm is the Grytczuk-Tropak algorithm.\\
\indent For every integer $v$ it is a consequence of Theorem
\ref{z} that there exists a minimal integer
$k$, $k_{\rm min}^{\epsilon}$, such that there exists a natural number $n$ with $a_n(k_{\rm min}^{\epsilon})=v$. 
Since ${\cal A}(k)={\cal C}(k)$ it follows that $k_{\rm min}^1=k_{\rm min}^{-1}$. Put 
$k_{\rm min}=k_{\rm min}^{-1}=k_{\rm min}^1$.
Grytczuk and Tropak \cite[Table 2.1]{GT} used their method to determine $k_{\rm min}$ for the integers in the
interval $[-9,\dots,10]$. Bosma 
\cite{Bosma} extended this to the range
$[-50,50]$ and we on our turn have extended the range from $[-70,\dots,70]$ (in which case
$k_{\rm min}\le 105$).

\section{Computation of $M(a_n(k))$ and $M(c_n(k))$}
Our starting point
is Lemma \ref{simpelzeg}. For each of the $p(k)$ summands we compute
the average (in Lemma \ref{halfweg}), which turns out to
be independent of $\epsilon$, and then find $M(a_n^{\epsilon}(k))$ by adding
these $p(k)$ averages. Note that of course $M(a_n(k))=M(a_n^{1}(k))$ and 
$M(c_n(k))=M(a_n^{-1}(k))$.
\begin{Lem}
\label{halfweg} Let $\lambda=(k_1^{n_{k_1}}\dots k_s^{n_{k_s}})$ be a partition with
$s\ge 1$, $k_1,\ldots,k_s$ distinct integers and $n_{k_1}\ge n_{k_2}\ge \dots \ge n_{k_s}\ge 1$.
If $n_{k_1}\ge 2$ we let $t$ be the largest integer $\le s$ for which $n_{k_t}\ge 2$, otherwise
we let $t=0$.
Let $L=[k_1,\dots,k_s],~G=(k_1,\dots,k_s)$ and $\epsilon=\pm 1$. We have
$$\lim_{x\rightarrow \infty}{1\over x}\sum_{n\le x}(-1)^{n_{k_1}+\cdots+n_{k_s}}\left({\epsilon\mu({n\over k_1})\atop n_{k_1}}\right)\cdots 
\left({\epsilon\mu({n\over k_s})\atop n_{k_s}}\right)={6\over \pi^2}{\epsilon_2(\lambda))
\over G\prod_{p|{L\over G}}(p+1)},$$
where
\begin{equation}
\label{diffie}
\epsilon_2(\lambda)=\epsilon_1(\lambda)\mu({L\over k_{t+1}})\cdots \mu({L\over k_s}),
\end{equation}
and
$$\epsilon_1(\lambda)=\cases{1 &if $n_{k_1}=1$, $s$ is even and $\mu(L/G)\ne 0$;\cr
\mu({L\over k_1})^{s-t}/2 &if $n_{k_1}\ge 2$ and $\mu(L/k_1)=\cdots=\mu(L/k_t)$ and $\mu(L/G)\ne 0$;\cr
0 & otherwise.}$$
\end{Lem}
Remark. A case by case analysis from
Lemma \ref{halfweg} on using (\ref{veertien}) shows that
$$2\epsilon_2(\lambda)=\prod_{j=1}^s (-1)^{n_{k_j}} \left({\mu({L\over k_j})\atop n_{k_j}}\right)
+  \prod_{j=1}^s (-1)^{n_{k_j}} \left({-\mu({L\over k_j})\atop n_{k_j}}\right).$$
\noindent {\it Proof of Lemma} \ref{halfweg}. 
Put 
$$S(x)=\sum_{n\le x}(-1)^{n_{k_1}+\cdots+n_{k_s}}\left({\epsilon\mu({n\over k_1})\atop n_{k_1}}\right)\cdots 
\left({\epsilon\mu({n\over k_s})\atop n_{k_s}}\right).$$
By (\ref{veertien}) for $n_{k_i}\ge 2$ the binomial coefficient $({\epsilon\mu(n/k_i)\atop n_{k_i}})$ 
is only non-zero 
if $\mu({n\over k_i})=-\epsilon$. 
Using (\ref{veertien}) it
follows that a necessary condition for the argument of $S(x)$ to be non-zero is that
$L|n$. Now write $n=mL$. Note that $\mu(mL/k_1)\cdots \mu(mL/k_s)=\mu(m)^s \mu(L/k_1)\cdots \mu(L/k_s)$
if $(m,L/k_j)=1$ for $1\le j\le s$ and equals zero otherwise. It is not difficult to show that
$[{L\over k_1},\dots,{L\over k_s}]={L\over G}$ and using this, that 
$\mu(L/k_i)\ne 0$ for $1\le i\le s$
iff $L/G$ is  squarefree. It follows
that if $\mu(L/G)=0$, then $S(x)=0$ and we are done, so next assume
that $\mu(L/G)\ne 0$. We infer that
$$S(x)=\sum_{m\le x/L,~(m,L/G)=1}(-1)^{n_{k_1}+\cdots+n_{k_s}}\left({\epsilon\mu(mL/k_1)\atop n_{k_1}}\right)\cdots 
\left({\epsilon\mu(mL/k_s)\atop n_{k_s}}\right).$$
Let us first consider the (easy) case where $n_{k_1}=1$. Then we obtain
$S(x)=(-1)^s\mu({L\over k_1})\cdots \mu({L\over k_s})\sum_{m\le x/L,~(m,L/G)=1}(\epsilon\mu(m))^s$.
If $s$ is odd, then by Lemma \ref{mobiustwee} it follows that $\lim_{x\rightarrow \infty}S(x)/x=0$
and we are done, so next assume that $s$ is even. Then we apply Lemma \ref{mobiuseen} and
obtain that
$$\lim_{x\rightarrow \infty}{S(x)\over x}={6 \mu({L\over k_1})\cdots \mu({L\over k_s})\over
\pi^2 L\prod_{p|{L\over G}}(1+{1\over p})}.$$
The assumption $\mu(L/G)\ne 0$ implies that  $L\prod_{p|L/G}(1+1/p)=G\prod_{p|L/G}(p+1)$.\\
Next we consider the case where $n_{k_1}\ge 2$. The corresponding binomial coefficient is
only non-zero if $\mu(mL/k_1)=-\epsilon$. Similarly, we must have
$\mu(mL/k_j)=-\epsilon$ for $1\le j\le t$. It follows that if it is not true that
$\mu(L/k_1)=\dots=\mu(L/k_t)$, then $S(x)=0$ and hence $\lim_{x\rightarrow \infty}S(x)/x=0$ 
as asserted, so assume that $\mu(L/k_1)=\dots=\mu(L/k_t)$. We have, on noting that
$\mu(mL/k_1)=-\epsilon$ and $(m,L/G)=1$ implies $-\epsilon\mu(m)=\mu(L/k_1)$,
\begin{eqnarray}
S(x)&=&\sum_{m\le x/L,~(m,L/G)=1\atop \mu(mL/k_1)=-\epsilon}(-\epsilon\mu(m))^{s-t}
\mu({L\over k_{t+1}})\cdots \mu({L\over k_s})\nonumber\\
&=&(\mu({L\over k_1}))^{s-t}\mu({L\over k_{t+1}})\cdots \mu({L\over k_s})
\sum_{m\le x/L,~(m,L/G)=1\atop \mu(m)=-\epsilon\mu(L/k_1)}1.\nonumber
\end{eqnarray}
On invoking Lemma \ref{mobiusdrie} the proof
is then completed. \qed

\begin{Thm}
\label{eenvoudig}
Let $\epsilon=\pm 1$. We have
$$M(a_n^{\epsilon}(k))={6\over \pi^2}\sum_{\lambda=(k_1^{n_{k_1}}\dots k_s^{n_{k_s}})\in {\cal P}(k)\atop
n_{k_1}\ge \dots \ge n_{k_s}\ge 1}{\epsilon_2(\lambda)\over G(\lambda)\prod_{p|{L(\lambda)\over G(\lambda)}}(p+1)},$$
where
$\epsilon_2(\lambda)$ is defined in {\rm (\ref{diffie})},
$L(\lambda)=[k_1,\dots,k_s]$ and $G(\lambda)=(k_1,\dots,k_s)$. Moreover, we have
$M(\alpha_n(k))=M(a_n(k))=M(c_n(k))$ and $M(\beta_n(k))=0$.
\end{Thm}
{\it Proof}. The first assertion follows from Lemma \ref{simpelzeg} together with
Lemma \ref{halfweg}. It is easy to see that $M(\alpha_n(k))$ and $M(\beta_n(k))$
exist. {}From 
$$M(\alpha_n(k))+M(\beta_n(k))=M(a_n(k))=M(c_n(k))=M(\alpha_n(k))-M(\beta_n(k)),$$
the second assertion then follows.\qed\\

\noindent In Table 3 we demonstrate Theorem \ref{eenvoudig} in case $k=4$.\\

\centerline{{\bf Table 3:} Computation of $e_4=\zeta(2)M(a_n(4))$}
\begin{center}
\begin{tabular}{|c|c|c|c|c|c|c|c|c|}
\hline
{\rm partition} & $\lambda$ & $n_{k_1}$ & $L(\lambda)$ & $G(\lambda)$ & $t$ & $s$ & $\epsilon_2(\lambda)$ & contribution to $e_4$\\
\hline
$4$  & $(4^1)$ & $1$ & $4$ & $4$ & $0$ & $1$ & $0$ & $0$\\
\hline
$3,1$  & $(3^11^1)$ & $1$ & $3$ & $1$ & $0$ & $2$ & $-1$ & $-1/4$\\
\hline
$2,2$  & $(2^2)$ & $2$ & $2$ & $2$ & $1$ & $1$ & $1/2$ & $+1/4$\\
\hline
$2,1,1$  & $(1^22^1)$ & $2$ & $2$ & $1$ & $1$ & $2$ & $-1/2$ & $-1/6$\\
\hline
$1,1,1,1$  & $(1^4)$ & $4$ & $1$ & $1$ & $1$ & $1$ & $1/2$ & $+1/2$\\
\hline
\end{tabular}
\end{center}
It is seen that $e_4=\zeta(2)M(a_n(4))=0-{1\over 4}+{1\over 4}-{1\over 6}+{1\over 2}={1\over 3}$.\\

\noindent The above formula suggests a connection with the group or representation
theory of the symmetric group $S_k$. The conjugacy classes in $S_k$ are in 1-1 correspondence
with the partitions of $k$. If $\lambda=(k_1^{n_{k_1}}\dots k_s^{n_{k_s}})$, then the order
of every element in the corresponding conjugacy class equals $L(\lambda)$. In particular,
$L(\lambda)\le g(k)$, where $g(k)$ denotes the maximum of all orders of elements in $S_k$.
It was shown by E. Landau in 1903 that $\log g(k)\sim \sqrt{k\log k}$ as $k$ tends to infinity
(for a nice account of this see \cite{Miller}, for recent results see
\cite{DNZ}), whereas by Stirling's theorem 
$\log k! \sim k\log k$. The average order of an element in $S_k$ is, not surprisingly, much
smaller than $g(k)$: if $\sigma$ is chosen at random from $S_k$, define
$Z=(\log {\rm order}(\sigma)-{1\over 2}\log^2 n)/(\log^{3/2}n/\sqrt{3})$, then the distribution
of $Z$ is known (see e.g. Nicolas \cite{Nicolas}) to converge to the standard normal 
distribution as $n\rightarrow \infty$.

\section{Average and value distribution}
\label{vvvier}
\noindent We give, using
Lemma \ref{naareindig2}, a simpler formula for $M(a_n(k))$ involving $a_n(k)$
for a finite set of $n$. 
\begin{Thm}
\label{vier}
Let $k\ge 1$ be fixed and $\epsilon=\pm 1$. Put $M_k=k\prod_{p\le k}p$. Then
$$M(a_n^{\epsilon}(k))={3\over \pi^2\prod_{p\le k}(1+{1\over p})}\sum_{d|M_k}{a_{d}^1(k)+a_{d}^{-1}(k)\over d}.$$
Furthermore, when $v\ne 0$,
$$\delta(a_n^{\epsilon}(k)=v)={3\over \pi^2\prod_{p\le k}(1+{1\over p})}\Big(\sum_{d|M_k\atop a_{d}^1(k)=v}{1\over d}+
\sum_{d|M_k\atop a_{d}^{-1}(k)=v}{1\over d}\Big).$$
\end{Thm} 
{\it Proof}. Let $r_1=\prod_{p\le k}p$. We have
\begin{eqnarray}
\sum_{n\le x}a_n^{\epsilon}(k)&=&\sum_{d|M_k}\sum_{dm\le x\atop (m,r_1)=1}(A_1(d)\mu(m)^2+\epsilon B_1(d)\mu(m))\nonumber\\
&=&\sum_{d|M_k}A_1(d)\sum_{m\le x/d\atop (m,r_1)=1}\mu(m)^2+o_k(x),\nonumber
\end{eqnarray}
where we used Lemma \ref{naareindig2} and Lemma \ref{mobiustwee}.
On invoking Lemma \ref{mobiuseen} we then obtain that 
$$\sum_{n\le x}a_n^{\epsilon}(k)={6x\over \pi^2\prod_{p\le k}(1+{1\over p})}\sum_{d|M_k}{A_1(d)\over d}+o_k(x).$$
On noting that
$A_1(d)=(a_{d}^1(k)+a_{d}^{-1}(k))/2$ (and $B_1(d)=(a_{d}^1(k)-a_{d}^{-1}(k))/2$, but this is not needed),
the first formula follows.\\
\indent As to the second identity we notice that, for $v\ne 0$, by Lemma \ref{naareindig2}
$$\sum_{n\le x,~a_n^{\epsilon}(k)=v}1=\sum_{d|M_k\atop a_d^1(k)=v}\sum_{md\le x,~\mu(m)=\epsilon\atop (m,r_1)=1}1
+\sum_{d|M_k\atop a_d^{-1}(k)=v}\sum_{md\le x,~\mu(m)=-\epsilon\atop (m,r_1)=1}1.$$
On invoking Lemma \ref{mobiusdrie}, the proof is then completed. \qed\\

\noindent Using identity (\ref{verdubbeling}) we arrive at the following corollary to this theorem:
\begin{cor}
\label{corvier}
Let $k\ge 3$ be fixed and odd and $M_k=k\prod_{p\le k}p$. Then
$$M(a_n^{\epsilon}(k))={1\over \pi^2\prod_{2<p\le k}(1+{1\over p})}\sum_{d|M_k/2}{a_{d}^1(k)+a_{d}^{-1}(k)\over d}.$$
Furthermore, when $v\ne 0$,
$$\delta(a_n^{\epsilon}(k)=v)={3\over \pi^2\prod_{p\le k}(1+{1\over p})}\Big(\sum_{d|M_k/2\atop a_{d}^1(k)
=v}{1\over d}+
\sum_{d|M_k/2\atop a_{d}^1(k)=-v}{1\over 2d}+
\sum_{d|M_k/2\atop a_{d}^{-1}(k)=v}{1\over d}+\sum_{d|M_k/2\atop a_{d}^{-1}(k)=-v}{1\over 2d}\Big).$$
\end{cor}
This result gives an alternative proof of the fact that, with $k\ge 3$ and odd, ${\cal A}(k)$ is symmetric.
Namely, it shows that for these $k$ we have $\delta(a_n^1(k)=v)>0$ iff $\delta(a_n^1(k)=-v)>0$.\\
\indent In case $k$ is prime, the divisor sum in the previous corollary can be further reduced.
\begin{Lem} 
\label{priempiet}
Let $k\ge 3$ be a fixed prime. Put $R_k=\prod_{2<p<k}p$. Then
$$M(a_n^{\epsilon}(k))={1\over \pi^2\prod_{2<p< k}(1+{1\over p})}\sum_{d|R_k}{a_{d}^1(k)+a_{d}^{-1}(k)\over d}.$$
\end{Lem}
{\it Proof}. We consider the formula given in the previous corollary. The divisors of
$M_k/2$ are either of the form $d$ with $d|kR_k$ or of the form $dk^2$ with $d|R_k$. For
the latter divisors $d$ we find, using (\ref{naarkwadraatvrijekern}) that
$a_{dk^2}^{\epsilon}(k)=a_{dk}^{\epsilon}(1)=-\epsilon\mu(dk)$ and hence
$\sum_{d|M_k/2}(a_d^1(k)+a_d^{-1}(k))/d=\sum_{d|kR_k}(a_d^1(k)+a_{d}^{-1}(k))/d$. Now
suppose that $d|R_k$. Using Lemma \ref{simpelzeg} we infer that
$$a_{dk}^{1}(k)+a_{dk}^{-1}(k)=a_d^1(k)+a_{d}^{-1}(k)-\mu(d)+\mu(d)=a_d^1(k)+a_{d}^{-1}(k).$$
Using this observation it follows that
$$\sum_{d|kR_k}{a_d^1(k)+a_{d}^{-1}(k)\over d}=(1+{1\over k})\sum_{d|R_k}{a_d^1(k)+a_{d}^{-1}(k)\over d},$$
whence the result follows. \qed
\begin{Lem}\label{ab} $~$\\
{\rm a)} We have $e_k2k\prod_{p\le k}(p+1)\in \mathbb Z$.\\
{\rm b)} If $k\ge 3$ is a prime, then $e_k2\prod_{p<k}(p+1)\in \mathbb Z$. 
\end{Lem}
{\it Proof}. a) An easy consequence of Theorem \ref{vier}.\\
b) An easy consequence of Lemma \ref{priempiet}. (It also follows from Theorem \ref{eenvoudigodd} below
on noting that $\epsilon_2(\lambda)=0$ in case $k|L(\lambda)$.)\qed\\

\noindent Numerically we observed that actually for $k\le 100$, we have $e_kk\prod_{p\le k}(p+1)\in \mathbb Z$. 
Using Lemma \ref{ab} it is seen 
that $e_kk\prod_{p\le k}(p+1)\in \mathbb Z$ if $k$ is an odd prime.\\
\indent Clearly $\delta(a_n(k)=0)=1-\sum_{v\ne 0}\delta(a_n(k)=v)$, where the latter
sum has only finitely many non-zero values and it is a finite computation to determine those $v$ for
which $a_n(k)=v$ for some $n$.   
For $1\le k\le 16$ the non-zero values of $\zeta(2)\delta(a_n(k)=v)$ are given
in Table 4 (except for $v=0$). By Table 4E we denote the extended version of Table 4, which
is Table 11 in \cite{HM}.\\

\centerline{{\bf Table 4:} Value  of $\zeta(2)\delta(a_n(k)=v)$}
\centerline{({\bf Table 4E}: Extended version of this table, see \cite[Table 11]{MH})}
\begin{center}
\begin{tabular}{|c|c|c|c|c|}
\hline
  & $v=-2$ & $v=-1$ & $v=1$ & $v=2$ \\
\hline
$k=1$  & $0$ & $1/2$ & $1/2$ & $0$ \\
\hline
$k=2$ & $0$ & $1/12$ & $7/12$ & $0$ \\
\hline
$k=3$ & $0$ & $5/24$ & $3/8$ & $0$ \\
\hline
$k=4$  & $0$ & $1/6$ & $1/2$ & $0$ \\ 
\hline
$k=5$  & $0$ & $13/80$ & $23/80$ & $0$ \\
\hline
$k=6$ & $0$ & $25/144$ & $67/144$ & $0$ \\
\hline
$k=7$ & $1/576$ & $577/2688$ & $731/2688$ & $1/1152$ \\
\hline
$k=8$  & $0$ & $1/8$ & $5/12$ & $0$ \\ 
\hline
$k=9$  & $0$ & $65/384$ & $347/1152$ & $0$ \\ 
\hline
$k=10$  & $0$ & $161/960$ & $347/960$ & $0$ \\ 
\hline
$k=11$  & $1/2304$ & $8299/50688$ & $11489/50688$ & $1/4608$ \\ 
\hline
$k=12$  & $0$ & $349/2304$ & $1009/2304$ & $0$ \\ 
\hline
$k=13$  & $43/48384$ & $219269/1257984$ & $277171/1257984$ & $43/96768$ \\ 
\hline
$k=14$  & $13/21504$ & $2395/21504$ & $2319/7168$ & $1/2304$ \\ 
\hline
$k=15$  & $13/32256$ & $1345/7168$ & $97247/322560$ & $13/64512$ \\ 
\hline
$k=16$  & $5/21504$ & $12149/64512$ & $1127/3072$ & $5/2688$ \\ 
\hline
\end{tabular}
\end{center}
\medskip
\noindent Let us now look at Theorem \ref{eenvoudig} and Theorem \ref{vier} from the viewpoint of
computational complexity. In Theorem \ref{eenvoudig} the sum has $p(k)$ terms and the estimate 
(\ref{harlit})
yields that $\log p(k)\sim \pi\sqrt{2k/3}$ as $k$ tends to
infinity. In Theorem \ref{vier} we sum over $t(k)$ terms where
$\log t(k)\sim \pi(k)\log 2 \sim k\log 2/\log k$. So Theorem \ref{eenvoudig} yields the
computational superior method. Theorem \ref{vier} is, however, much more easily implemented.
Using Lemma \ref{simpelzeg} in combination with Theorem \ref{vier} and 
Lemma \ref{omi} below, an alternative
proof of Theorem \ref{eenvoudig} is obtained. If one starts with Lemma \ref{priempiet}
and invokes Lemmas \ref{simpelzeg} and \ref{omi}, one obtains a sum over partitions of $k$, where now
only odd integers are allowed to occur in the partition. This yields a result superior in complexity
to that provided by Theorem \ref{eenvoudig}, since for $p_{\rm odd}(k)$ the number of partitions of $m$ into
odd parts we have $\log p_{\rm odd}(k)\sim \pi\sqrt{k/3}$ (see 
e.g. Bringmann \cite{bring}), whereas $\log p(k)\sim \pi\sqrt{2k/3}$.
\begin{Thm}
\label{eenvoudigodd}
Let $\epsilon=\pm 1$ and $k\ge 3$ be a prime. With the notation from Theorem {\rm \ref{eenvoudig}} we have
$$M(a_n^{\epsilon}(k))={2\over \pi^2}\sum_{\lambda=(k_1^{n_{k_1}}\dots k_s^{n_{k_s}})\in {\cal P}(k)\atop
n_{k_1}\ge \dots n_{k_s}\ge 1}{\epsilon_2(\lambda)\over \prod_{p|L(\lambda)}(p+1)},$$
where
the sum is over all partitions of $k$ into only odd parts having at least one number repeated more
than once (i.e. $n_{k_1}\ge 2$).
\end{Thm}
The restriction that $n_{k_1}\ge 2$ does not yield an extra asymptotical improvement: using a result of
Hagis \cite{Hagis}, one sees that with $p_1(k)$ the number of partitions of $k$ into only odd parts 
having at least one number repeated more
than once (i.e. $n_{k_1}\ge 2$), we have $p_1(k)\sim p_{\rm odd}(k)$.\\
\indent Indeed, all results involving $\sum_{d|r}a_d^{\epsilon}(k)/d$ can be turned into a M\"oller type
of result involving partitions of $k$ on invoking the following lemma and Lemma \ref{simpelzeg}.
\begin{Lem}
\label{omi}
Let $L=[k_1,\ldots,k_s]$ and $G=(k_1,\ldots,k_s)$. The sum
$$\sum_{d|r}{(-1)^{n_{k_1}+\cdots+n_{k_s}}\over d}\left({\epsilon\mu({d\over k_1})\atop n_{k_1}}\right)\cdots 
\left({\epsilon\mu({d\over k_s})\atop n_{k_s}}\right)$$
equals $${(-\epsilon)^s\over L}\mu({L\over k_{t+1}})\cdots \mu({L\over k_s})\prod_{p|{r\over L},~p\nmid {L\over G}}\Big(1+{(-1)^s\over
p}\Big)$$
if $n_{k_1}=1$, $\mu({L/G})\ne 0$ and $r|L$,
it equals 
$${1\over 2L}\mu({L\over k_1})^{s-t}
\mu({L\over k_{t+1}})\cdots \mu({L\over k_s})
\Big[
\prod_{p|{r\over L},~p\nmid {L\over G}}(1+{1\over p})-\epsilon \mu({L\over k_1})
\prod_{p|{r\over L},~p\nmid {L\over G}}(1-{1\over p})\Big]$$
if $n_{k_1}\ge 2$, $\mu(L/k_1)=\cdots=\mu(L/k_t)\ne 0$ and $r|L$, and zero in all other cases.
\end{Lem}
\begin{cor} 
\label{blubber}
We have
$${1\over 2}\sum_{d|r}{(-1)^{n_{k_1}+\cdots+n_{k_s}}\over d}\Big(\left({-\mu({d\over k_1})\atop n_{k_1}}\right)\cdots 
\left({-\mu({d\over k_s})\atop n_{k_s}}\right)+
\left({\mu({d\over k_1})\atop n_{k_1}}\right)\cdots 
\left({\mu({d\over k_s})\atop n_{k_s}}\right)\Big)$$
$$
=\cases{{\epsilon_2(\lambda)\over L(\lambda)}\prod_{p|{r\over L(\lambda)},~p\nmid {L(\lambda)\over G(\lambda)}}(1+{1\over p}) & if
$L(\lambda)|r$;\cr
0 & otherwise.}$$
In particular, if  $r$ is squarefree, then the sum equals
$$
=\cases{\epsilon_2(\lambda)\prod_{p|r}(1+{1\over p})\prod_{p|L(\lambda)}(p+1)^{-1} & if $L(\lambda)|r$;\cr
0 & otherwise.}$$
\end{cor}
{\it Proof of Lemma} \ref{omi}. The proof is similar to 
that of Lemma \ref{halfweg}. Indeed, if $S(x)$ now denotes the new sum under consideration, then
the proof proceeds as that of Lemma \ref{halfweg}. Instead of the sum
$$\sum_{m\le x/L,~(a,L/G)=1}(\epsilon \mu(m))^s~{\rm ~we~have~}{1\over L}\sum_{mL|r,~(m,L/G)=1}
{(\epsilon \mu(m))^s\over m}.$$
Instead of the sum
$$\sum_{m\le x/L,~(m,L/G)=1\atop \mu(m)=-\epsilon \mu(L/k_1)}1
~{\rm ~we~have~}{1\over L}\sum_{mL|r,~(m,L/G)=1\atop \mu(m)=-\epsilon \mu(L/k_1)}{1\over m}.$$
On relating these sums to Euler products, the result follows. \qed\\

\noindent {\it Alternative proof of Theorem} \ref{eenvoudig}. On combining Theorem \ref{vier}, Lemma 
\ref{simpelzeg} and
Corollary \ref{blubber} (with $r=M_k$), we find that
$$M(a_n^{\epsilon}(k))= 
{6\over \pi^2\prod_{p\le k}(1+{1\over p})}\sum_{\lambda=(k_1^{n_{k_1}}\dots k_s^{n_{k_s}})\in {\cal P}(k)\atop 
L(\lambda)|M_k}{\epsilon_2(\lambda)\over L(\lambda)}\prod_{p|{M_k\over L(\lambda)},~p\nmid {L(\lambda)\over 
G(\lambda)}}(1+{1\over p}).$$
Now 
$$\prod_{p|{M_k\over L(\lambda)},~p\nmid {L(\lambda)\over 
G(\lambda)}}(1+{1\over p})={\prod_{p|{M_k\over G(\lambda)}}(1+{1\over p})\over 
\prod_{p|{L(\lambda)\over G(\lambda)}}(1+{1\over p})}={\prod_{p\le k}(1+{1\over p})\over \prod_{p|{L(\lambda)\over
G(\lambda)}}(1+{1\over p})},$$
where we used that
$$\prod_{p|u,~p\nmid v}(1+{1\over p})=\prod_{p|uv}(1+{1\over p})\prod_{p|v}(1+{1\over p})^{-1}$$
and $G(\lambda)|k$.  It thus follows that
\begin{equation}
\label{frankie}
M(a_n^{\epsilon}(k))= 
{6\over \pi^2}\sum_{\lambda=(k_1^{n_{k_1}}\dots k_s^{n_{k_s}})\in {\cal P}(k)\atop 
L(\lambda)|M_k}{\epsilon_2(\lambda)\over L(\lambda)}\prod_{p|{L(\lambda)\over G(\lambda)}}(1+{1\over p})^{-1}.
\end{equation}
If $\epsilon_2(\lambda)\ne 0$, then $L(\lambda)/G(\lambda)$ is squarefree and hence
$$G(\lambda){L(\lambda)\over G(\lambda)}\prod_{p|{L(\lambda)\over G(\lambda)}}(1+{1\over
p})=G(\lambda)\prod_{p|{L(\lambda)\over G(\lambda)}}(p+1)$$
and furthermore $L(\lambda)|M_k$. These two observations in combination with (\ref{frankie}) complete
the alternative proof of Theorem \ref{eenvoudig}. \qed\\

\noindent {\it Proof of Theorem} \ref{eenvoudigodd}. On combining 
Lemma \ref{priempiet}, Lemma \ref{simpelzeg} and Corollary \ref{blubber} 
(with $r=R_k$ a squarefree number), we find that 
\begin{eqnarray}
M(a_n^{\epsilon}(k))&=&{2\over \pi^2\prod_{2<p<k}(1+{1\over p})}\sum_{\lambda=(k_1^{n_{k_1}}\dots k_s^{n_{k_s}})\in {\cal P}(k)\atop
L(\lambda)|R_k}\epsilon_2(\lambda){\prod_{p|R_k}(1+{1\over p})\over 
\prod_{p|L(\lambda)}(p+1)};\nonumber\cr
&=&{2\over \pi^2}\sum_{\lambda=(k_1^{n_{k_1}}\dots k_s^{n_{k_s}})\in {\cal P}(k)\atop
L(\lambda)|R_k}{\epsilon_2(\lambda)\over 
\prod_{p|L(\lambda)}(p+1)}.\nonumber
\end{eqnarray}
Since $G(\lambda)|k$ and by assumption $k$ is an odd prime, either
$G(\lambda)=1$ or $G(\lambda)=k$. The latter case only occurs if $\lambda=(k^1)$ in which case
$\epsilon_2(\lambda)=0$, hence we may assume that $G(\lambda)=1$. Let us assume that
$\epsilon_2(\lambda)\ne 0$ and so $\mu(L(\lambda)/G(\lambda))\ne 0$ and so $L(\lambda)$ must
be squarefree. Thus each part $k_i$ of such a partition is squarefree and since $k_i\le k$ it
follows that $k_i|R_k$ iff $k_i$ is odd.
We infer that if $2\nmid L(\lambda)$, then $L(\lambda)|R_k$. If $L(\lambda)$ is
even, then $L(\lambda)\nmid R_k$. Thus the sum over all partitions with $L(\lambda)|R_k$, can
be restricted to those partitions consisting of only odd parts. If $n_{k_1}=1$, then 
the partition consists of distinct odd parts and so the number of parts $s$ must be odd, as
by assumption $k$ is odd, and hence $\epsilon_2(\lambda)=0$ in this case. Thus we can further
resctrict our partition sum to the partitions into odd parts only having $n_{k_1}\ge 2$. \qed

\section{Some observations related to Table 4}
In this section we make some observations regarding Table 4 (and Table 4E) and prove some
results inspired by these observations.\\
\indent For $k$ is even numerical results suggest that often ${\cal A}(k)$ is
not symmetric, whereas we have shown (Lemma \ref{symmie}) that for $k$ is odd it is always
symmetric. For $k\le 100$ it is mostly true that if $v\in {\cal A}(k)$ and $v$ is negative, then
$-v\in {\cal A}(k)$. This leads to the question as to whether perhaps 
$A_{+}(k)>A_{-}(k)$ for all $k$ sufficiently large, with $A_{+}(k)=\max {\cal A}(k)$ and $A_{-}(k)=-\min {\cal A}(k)$. 
It has been shown by Bachman \cite{gennady} that (\ref{gennady}) also holds true if we replace $A(k)$ by
$A_{+}(k)$, or $A_{-}(k)$. However, this result is not strong enough to decide on the above question.\\
\indent An other observation that can be made is that for $k\le 100$ it is
true that ${\cal A}(k)$ is convex, that is consists of consecutive integers, i.e. if $v_0<v_1$ are
in ${\cal A}(k)$, then so are all integers between $v_0$ and $v_1$.\\
\indent Let us define  ${\cal A}_j(k)=\{a_n(k):n\equiv j({\rm mod~}2)\}$, for $0\le j\le 1$.
\begin{Lem}$~$\\
{\rm 1)} We have $0\in {\cal A}_j(k)$.\\
{\rm 2)} If $k$ is even, then ${\cal A}_1(k)\subseteq {\cal A}_0(k)={\cal A}(k)$.\\ 
{\rm 3)} If $k$ is odd, then $v\in {\cal A}_1(k)$ iff $-v \in {\cal A}_0(k)$. \\
\end{Lem}
{\it Proof}. 1) Consider any integer $n_j$ such that $n_j/\gamma(n_j)>k$ and $n_j\equiv 
j({\rm mod~}2)$. Then, by part 1 of Lemma \ref{tralala}, we have $a_{n_j}(k)=0$ and hence
$0\in {\cal A}_j(k)$.\\
2) If $v\in {\cal A}_1(k)$, then $v=a_d(k)$ for some odd integer $d$. Then, by (\ref{verdubbeling})
we have $a_{2d}(k)=(-1)^ka_{d}(k)=v$ and
hence $v\in {\cal A}_0(k)$. We have ${\cal A}(k)={\cal A}_0(k)\cup {\cal A}_1(k)={\cal A}_0(k)$,
since ${\cal A}_1(k)$ is included in ${\cal A}_0(k)$.\\
3) Proceding as in part 2 we infer that if $v\in {\cal A}_1(k)$, then $-v\in {\cal A}_0(k)$. For the
reverse implication we make use of Lemma \ref{valueset}.\qed\\

\noindent Inspection of Table 4E shows that for odd integers $k$ with $A(k)\ge 2$ often 
$\delta(a_n(k)=A(k))$ and $\delta(a_n(k)=-A(k))$ differ by
a factor two or a factor less than two. Regarding this situation we have the following
result:
\begin{Lem}
\label{factortwo}
Let $k\ge 3$ be odd and $v\ne 0$. We have 
$${1\over 2}\delta(a_n^{\epsilon}(k)=v)\le \delta(a_n^{\epsilon}(k)=-v)\le 2\delta(a_n^{\epsilon}(k)=v).$$
Furthermore, we have
$$2\delta(a_n^{\epsilon}(k)=v)=\delta(a_n^{\epsilon}(k)=-v){\rm ~iff~}v\not\in {\cal A}_1(k) 
~({\rm that~is~iff~} -v\not\in {\cal A}_0(k)).$$
\end{Lem}
{\it Proof}. We write $w_k={3\over \pi^2}\prod_{p\le k}(1+{1\over p})^{-1}$ and (with $\epsilon_1=\pm 1$)
$$\alpha_{\epsilon_1}=\sum_{d|M_k/2\atop a_d^{1}(k)=\epsilon_1 v}{1\over d}+
\sum_{d|M_k/2\atop a_d^{-1}(k)=\epsilon_1 v}{1\over d}.$$
Note that $\alpha_{\epsilon_1}\ge 0$. Then, by Corollary \ref{corvier}, we have
$$\delta(a_n^{\epsilon}(k)=v)=w_k(\alpha_1+\alpha_{-1}/2) {\rm ~and~} 
\delta(a_n^{\epsilon}(k)=-v)=w_k(\alpha_{-1}+\alpha_{1}/2).$$ 
The first part of the assertion follows on comparing these two formulae.
As to the second assertion, the latter formulae imply that it is enough to prove that
$\alpha_1=0$ iff $v\not\in {\cal A}_1(k)$. A minor variation of
the proof of Lemma \ref{brovo} shows that $\{-1,0,1\}\subseteq \{a_n^1(k)|n>1,~2\nmid n\}$
and hence ${\cal A}_1(k)=\{a_n^1(k):2\nmid n\}$.   
A minor modification of the proof of Lemma \ref{valueset} shows that 
$\{a_n^1(k):2\nmid n\}=\{0,a_d^1(k),a_d^{-1}(k):d|M_k/2\}$.  
On noting that $M_k/2$ is odd, we infer that if
$v\not\in {\cal A}_1(k)$, then clearly $\alpha_1=0$. On the other hand, if
$\alpha_1=0$, then $v$ is not in $\{0,a_d^1(k),a_d^{-1}(k):d|M_k/2\}={\cal A}_1(k)$, completing
the proof. \qed \\

\centerline{{\bf Table 5:} Set theoretic difference ${\cal A}(k)\backslash {\cal A}_0(k)$ in 
case ${\cal A}(k)\ne {\cal A}_0(k)$ and $k\le 53$}
\begin{center}
\begin{tabular}{|c|c|c|c|c|}
\hline
$k=7$    & $\{-2\}$ & $k=11$ & $\{-2\}$ \\
\hline
$k=13$ & $\{-2\}$ & $k=15$ & $\{-2\}$ \\
\hline
$k=17$ & $\{-3\}$ & $k=19$ & $\{-3\}$ \\
\hline
$k=21$  & $\{-3\}$ & $k=23$ & $\{-4,-3\}$ \\
\hline
$k=25$  & $\{-3\}$ & $k=31$ & $\{-4\}$ \\
\hline
$k=35$  & $\{5\}$ & $k=37$ & $\{5\}$ \\
\hline
$k=39$  & $\{5,6\}$ & $k=43$ & $\{-7\}$ \\
\hline
$k=45$  & $\{-7\}$ & $k=47$ & $\{-9,-8\}$ \\
\hline
$k=51$  & $\{8\}$ & $k=53$ & $\{9,10,11,12,13\}$ \\
\hline
\end{tabular}
\end{center}
\medskip
\noindent Example. Inspection of Table 4 shows that 
$\delta(a_n(7)=-2)=2\delta(a_n(7)=2)$. It thus follows by Lemma \ref{factortwo} that there is no even integer $n$ for which
$a_n(7)=-2$ (whereas $a_{105}(7)=-2$). Further examples can be derived from Table 5.\\

\noindent For $k\le 100$ it turns out ${\cal A}(k)\backslash {\cal A}_0(k)$ is always {\it convex}, i.e.
consists of consecutive integers. 
For part 2 of Lemma \ref{factortwo} to be of some mathematical value we would hope that
  infinitely often ${\cal A}_0(k)$ is strictly contained in ${\cal A}(k)$. Note
  that by Theorem \ref{z} we have $\{{\cal A}_0(k):k\ge 1\}=\{{\cal A}(k):k\ge 1\}=\mathbb Z$. \\
\indent Suppose that $\delta(a_n(k)=v)>0$. Then the quotient 
$${\delta(a_n(k)=-v)\over \delta(a_n(k)=v)}$$
does not exceed 2 in case $k$ is odd by Lemma \ref{factortwo}. Inspection of Table 4E 
suggests that given {\it any} real number $r$, we can find
a $v>0$ and even $k$ such that the latter quotient exceeds $r$.

\section{Some variations}
Using the same methods we can easily determine e.g. $M(\mu(n)^2a_n(k))$, i.e. the average of $a_n(k)$ over
all squarefree integers $n$. We obtain the following result.
\begin{Thm}
\label{vierkwadraatvrij}
Let $k\ge 1$ be fixed and put $Q_k=\prod_{p\le k}p$. Then
$$M(\mu(n)a_n^{\epsilon}(k))={3\epsilon\over \pi^2\prod_{p\le k}(1+{1\over p})}\sum_{d|Q_k}{\mu(d)(a_{d}^1(k)-
a_{d}^{-1}(k))\over d},$$
and
$$M(\mu(n)^2a_n^{\epsilon}(k))={3\over \pi^2\prod_{p\le k}(1+{1\over p})}\sum_{d|Q_k}{a_{d}^1(k)+a_{d}^{-1}(k)\over d}.$$
\end{Thm} 
This result implies that $M(\mu(n)a_n(k))=-M(\mu(n)c_n(k))$.
\begin{Thm}
\label{eenvoudig2}
Let $\epsilon=\pm 1$. We have
$$M(\mu(n)^2a_n^{\epsilon}(k))={6\over \pi^2}\sum_{\lambda=(k_1^{n_{k_1}}\dots k_s^{n_{k_s}})\in {\cal P}(k)}
{\epsilon_2(\lambda)\mu(L(\lambda))^2\over \prod_{p|L(\lambda)}(p+1)}.$$
\end{Thm}
{\it Proof}. A simple variation of the proof of Theorem \ref{eenvoudig}.\qed\\

\noindent Put $$f_k=\zeta(2)M(\mu(n)a_n(k)){\rm ~and~}g_k=\zeta(2)M(\mu(n)^2a_n(k)).$$
Note that
$$g_k=\lim_{x\rightarrow \infty}{\sum_{n\le x}\mu(n)^2a_n(k)\over \sum_{n\le x}\mu(n)^2}.$$
$~$
\vfil\eject
\centerline{{\bf Table 6:} Scaled averages, $f_k=\zeta(2)M(\mu(n)a_n(k))$ and $g_k=\zeta(2)M(\mu(n)^2a_n(k))$}
\begin{center}
\begin{tabular}{|c|c|c|c|c|c|c|c|c|c|c|c|}
\hline
$k$  & $1$ & $2$ & $3$ & $4$ & $5$ & $6$ & $7$ & $8$ & $9$ & $10$ \\
\hline
\T \B $f_k$  & $-1$ & $-{1\over 6}$ & $-{1\over 4}$ & $-{1\over 6}$ & $-{5\over 24}$ & 
$-{1\over 12}$ & $-{7\over 24}$ & $-{7\over 72}$ & $-{7\over 48}$ & $-{7\over 72}$\\
\hline
\T \B $g_k$  & $0$ & ${1\over 2}$ & ${1\over 6}$ & ${1\over 4}$ & ${1\over 8}$ & 
${1\over 3}$ & ${1\over 18}$ & ${5\over 24}$ & ${17\over 144}$ & ${23\over 96}$\\
\hline
$k$  & $11$ & $12$ & $13$ & $14$ & $15$ & $16$ & $17$ & $18$ & $19$ & $20$\\
\hline
\T \B $f_k$  & $-{25\over 96}$ & $-{31\over 576}$ & $-{11\over 42}$ & $-{1\over 16}$ & $-{11\over 84}$ &
$-{8\over 63}$ & $-{491\over 2688}$ & $-{613\over 12096}$ & $-{2371\over 10080}$ & $-{173\over 4032}$  \\
\hline
\T \B $g_k$  & ${1\over 16}$ & ${59\over 192}$ & ${13\over 288}$ & ${155\over 672}$ & ${145\over 1344}$ & 
${425\over 4032}$ & ${667\over 8064}$ & ${523\over 2016}$ & ${55\over 1344}$ & ${101\over 480}$\\
\hline
\end{tabular}
\end{center}
\medskip
\begin{Lem}$~$\\
{\rm 1)} If $k$ is a prime, then $g_k=e_k$.\\
{\rm 2)} Let $k=2q$ with $q$ an odd prime. Then 
$$g_{2q}=e_{2q}+{e_q\over 2}-{1\over 2q(q+1)}.$$
\end{Lem}
{\it Proof}. 1) Reasoning as in the beginning of the proof of Lemma \ref{priempiet} we
infer that in case $k>2$ is prime, we have
$$M(a_n(k))={3\over \pi^2\prod_{p\le k}(1+{1\over p})}\sum_{d|Q_k}{a_{d}^1(k)+a_{d}^{-1}(k)\over d}
=M(\mu(n)^2a_n(k)).$$
Since $e_2=g_2={1\over 2}$, the proof is completed.\\
2) Note that $G(\lambda)|2q$. We consider the contribution of the $\lambda\in {\cal P}(q)$ with
$G(\lambda)=1,2,q$ or $2q$ separately. Denote these by, respectively, 
$\Sigma_1,\Sigma_2,\Sigma_{q}$ and $\Sigma_{2q}$. In case $G(\lambda)=2$, we let 
$k'_i=k_i/2$ for $1\le i\le s$ and let $\lambda'=({k'_1}^{n_{k'_1}}\dots {k'_s}^{n_{k'_s}})$ be
the associated partition of $q$. We have $G(\lambda')=1$, $L(\lambda)=2L(\lambda')$ and
$\epsilon_2(\lambda)=\epsilon_2(\lambda')$. Thus
$$\Sigma_2={3\over \pi^2}\sum_{\lambda'=({k'_1}^{n_{k'_1}}\dots {k'_s}^{n_{k'_s}})\in {\cal P}(q)}
{\epsilon_2(\lambda')\over \prod_{p|L(\lambda')}(p+1)}={1\over 2}M(a_n(q)),$$
by Theorem \ref{eenvoudig}. It is easily seen that $\Sigma_3={1\over 2q\zeta(2)}$ and $\Sigma_{2q}=0$.
On putting everything together we obtain
\begin{equation}
\label{padua}
M(a_n(2q))=\Sigma_1+{1\over 2}M(a_n(q))+{1\over 2q\zeta(2)}.
\end{equation}
Likewise we write $M(\mu(n)^2a_n(2q))=\Sigma'_1+\Sigma'_2+\Sigma'_q+\Sigma'_{2q}$. 
It is easily seen that
\begin{equation}
\label{paris}
M(\mu(n)^2a_n(2q))=\Sigma_1+\Sigma'_2+{1\over 2(q+1)\zeta(2)}.
\end{equation}
We write
$\Sigma'_2=\Sigma'_{2,1}+\Sigma'_{2,2}$, with $\lambda$ contributing to $\Sigma'_{2,1}$ if
it contributes to $\Sigma'_2$ and $4\nmid L(\lambda)$, and $\lambda$ contributing to $\Sigma'_{2,2}$ if
it contributes to $\Sigma'_2$ and $4|L(\lambda)$. 
As before to a $\lambda\in {\cal P}(2q)$ with $G(\lambda)=2$, we associate a partition $\lambda'\in {\cal P}(q)$. Note
that $G(\lambda')=1$, $L(\lambda)=2L(\lambda')$ and $\epsilon_2(\lambda')=\epsilon_2(\lambda)$.
On invoking Theorem \ref{eenvoudig2} we find that
\begin{eqnarray}
\Sigma'_{2,1}&=&{6\over \pi^2}\sum_{{\lambda=(k_1^{n_{k_1}}\dots k_s^{n_{k_s}})\in {\cal P}(2q)\atop G(\lambda)=2,~4\nmid
L(\lambda)}}{\epsilon_2(\lambda)\mu(L(\lambda))^2\over \prod_{p|L(\lambda)}(p+1)}\nonumber\cr
&=&{2\over \pi^2}\sum_{{\lambda'=({k'_1}^{n_{k'_1}}\dots {k'_s}^{n_{k'_s}})\in {\cal P}(q)\atop G(\lambda')=2,~2\nmid
L(\lambda')}}{\epsilon_2(\lambda')\over \prod_{p|L(\lambda')}(p+1)}=M(a_n(q)),\nonumber\cr
\end{eqnarray}
by Theorem \ref{eenvoudigodd}. Since obviously $\Sigma'_{2,2}=0$, it follows from (\ref{paris}) that
$$M(\mu(n)^2a_n(2q))=\Sigma_1+M(a_n(q))+{1\over 2(q+1)\zeta(2)}.$$
The result follows on equating $\Sigma_1$ coming from the latter equality with $\Sigma_1$ coming
from (\ref{padua}). \qed\\

\noindent Remark. An alternative proof of part 1 of the latter lemma is obtained on noting that if $k$ is an odd
prime, then
$${\epsilon_2(\lambda)\mu(L(\lambda))^2\over \prod_{p|L(\lambda)}(p+1)}={\epsilon_2(\lambda)\over
G(\lambda)\prod_{p|{L(\lambda)\over G(\lambda)}}(p+1)},$$
and invoking Theorem \ref{eenvoudig2} and Theorem \ref{eenvoudig}.

\section{Open problems}
For the convenience of the reader we have collected below the open problems arising in this paper.\\
({\bf P1}) Is it true that ${\cal A}(k)$ is convex ?\\
({\bf P2}) Is it true that ${\cal A}(k)\backslash {\cal A}_0(k)$ is convex ?\\
({\bf P3}) Is it true that $\epsilon M(\mu(n)a_n^{\epsilon}(k))<0$ ?\\
({\bf P4}) Is it true that $M(\mu(n)^2a_n^{\epsilon}(k))>0$ for $k\ge 2$ ?\\
({\bf P5}) Is M\"ollers conjecture that $0\le e(k)\le 1/2$ true ?\\
({\bf P6}) Is $e_kk\prod_{p\le k}(p+1)$ always an integer ? (Certainly true if $k$ is an odd prime.)\\
({\bf P7}) What can one say about the behaviour of $e(k)$ as $k$ gets large, or $k$ has many
distinct prime factors ?\\
({\bf P8}) What is the smallest integer $k_0$ such that $A(k_0)>k_0$ ? M\"oller \cite[(10)]{HM} has shown that
$k_0\le 1820$. Our computations show that $k_0>105$.\\
({\bf P9}) Find effective estimates for $A(k)$.\\
({\bf P10}) Is it true that infinitely often $A(k)>A_0(k)$ ?\\
({\bf P11}) Is it true that $A_{+}(k)>A_{-}(k)$ for all $k$ sufficiently large ?\\
({\bf P12}) Given any real number $r$, can we find $k$ and $v$ such that $\delta(a_n(k)=v)\ne 0$ 
and $\delta(a_n(k)=-v)>r\delta(a_n(k)=v)$ ?\\
({\bf P13})  Determine $\{a_d^{\epsilon}(k)~|~d|M_k\}$, cf. p. 10.\\

\vfil\eject
\noindent {\bf Acknowledgement}. This paper is partly based on the M.Sc.~thesis of
Huib Hommersom (University of Amsterdam). This was a `literature M.Sc. thesis', i.e. the student was assigned
to get an overview of the literature in a specific topic and report on this. A
version of this thesis including research contributions (all by Moree) was posted on the
arXiv in 2003 \cite{MH}. The extensive numerical work in \cite{MH} was carried out by Y. Gallot
and partially checked using Maple by Moree.\\
\indent The present paper was written whilst the second author was working at the
Max-Planck-Institute for Mathematics in Bonn. He likes to thank both this institute and the University
of Amsterdam for their pleasant and inspiring research atmospheres.

\medskip\noindent {\footnotesize 12 bis rue Perrey,\\
31400 Toulouse, France.\\
e-mail: {\tt galloty@orange.fr}}\\

\medskip\noindent {\footnotesize Max-Planck-Institut f\"ur Mathematik,\\
Vivatsgasse 7, D-53111 Bonn, Germany.\\
e-mail: {\tt moree@mpim-bonn.mpg.de}}\\

\medskip\noindent {\footnotesize Bestevaerstraat 46, 1056 HP Amsterdam, The Netherlands\\
e-mail: {\tt hjhom48@hotmail.com}}
\end{document}